%% file: Comments_on_A_Framework.tex
\documentclass[
 onecolumn,
 11pt, 
 draftcls,
journal
]{IEEEtran}

\usepackage{array}
\usepackage{epsfig}
\usepackage{cite}
\usepackage{amsfonts}
\usepackage{amsmath}
\usepackage{amssymb}
\usepackage{amsthm}
\usepackage{amsxtra}
\usepackage{varioref}
\usepackage{multicol}
\usepackage{float}
\usepackage{rotate}
\usepackage{color}
\usepackage{enumerate}
\usepackage[active]{srcltx} % this should allow forward search
\usepackage{mathrsfs}
\usepackage{color}
\usepackage{rotating}
\usepackage{times}
\usepackage[english]{babel}
\usepackage[utf8]{inputenc}
\usepackage{fontenc}
\usepackage{Milancommands}
\usepackage{tikz}
  \usetikzlibrary {calc,positioning,graphs,quotes,arrows.meta,math}
\tikzset{
>={Latex[width=3pt,length=4pt]},
 bloq/.style={rectangle, rounded corners=.7mm, inner ysep=1.3pt, minimum width=5mm, draw=black, font=\itshape},
 sm/.style={circle, draw=black,inner sep=.7pt},
 cov/.style 2 args={to path={-- ++(#1,0) node[pos=.5,above] {#2}|- (\tikztotarget)}},
 hv/.style 2 args={to path={-| (\tikztotarget) node[left,pos=#1]{#2}}},
 vh/.style={to path={|- (\tikztotarget)}}
 }
\usepackage[colorlinks=true,linkcolor=black,citecolor=black]{hyperref}
\newcommand{\enc}{\Esp}
\newcommand{\dec}{\Dsp}
\newcommand{\EC}{\text{EC}}

\newcommand{\pla}{\mathscr{F}}
\renewcommand{\rvar}{r}
\renewcommand{\rvap}{p}
\renewcommand{\rvaq}{q}
\renewcommand{\rvas}{s}
\renewcommand{\rvae}{e}
\renewcommand{\rvax}{x}
\renewcommand{\rvay}{y}
\renewcommand{\rvau}{u}
\renewcommand{\rvaz}{z}
\newtheorem{thm}{Theorem}
\newtheorem{defn}{Definition}

\newtheorem{rem}{Remark}

\title{Comments on\\ ``A Framework for Control System Design Subject to Average Data-Rate Constraints''} 

 \author{ Milan S. Derpich$^{\diamond}$
 \thanks{$^{\diamond}$Department of  Electronic Engineering,
 Universidad T\'ecnica Federico Santa Mar\'ia, Valpara\'iso, Chile.
milan.derpich@usm.cl.}
and 
Jan {\O}stergaard$^{\dagger}$
\thanks{$^{\dagger}$Department of Electronic Systems, Aalborg University, Fredrik Bajers Vej 7, DK-9220, Aalborg, Denmark; janoe@ieee.org.}}

\begin{document}
\maketitle
\begin{abstract}
Theorem~4.1 in the 2011 paper ``A Framework for Control System Design Subject to Average Data-Rate Constraints''
allows one to lower bound average operational data rates in feedback loops (including the situation in which encoder and decoder have side information). 
Unfortunately, its proof is invalid.

In this note we first state the theorem and explain why its proof is flawed, and then  provide a correct proof under weaker assumptions.

\end{abstract}

\section{Introduction}
The paper~\cite{silder11}  provides a set of theorethical tools to design networked stochastic control problems by using linear time invariant systems and entropy-coded dithered quantizers.
In this paradigm, the average data rate is precisely linked to the signal-to-noise ratio in an  identical system where the encoder and decoder are replaced by an additive white Gaussian noise channel.
In this way, the design task can be carried out in a stochastic control arena.

 In this note we reveal that the proof of~\cite[Theorem 4.1]{silder11}, a key result of that paper, is flawed.
 Then we use a recent result by the authors~\cite[Theorem~2]{derost21a} 
 not only to derive a simple (and correct) proof of~\cite[Theorem 4.1]{silder11}, but also to generalize it.

\section{Notation and Preliminaries}
If $x$ denotes a random sequence, then 
$x=x(0),x(1),\ldots$
and
$x^{k}\eq x(0),x(1),\ldots,x(k)$.

Denote the marginal probability distributions of two random variables $\rvax$ and $\rvay$ by
$\mu$,
$\nu$, respectively.
Define the product measure $\pi\eq \mu\times\nu$.
The mutual information between $\rvax$ and $\rvay$ is defined as
 \begin{align}
  I(\rvax;\rvay)\eq  \int \log\left( \frac{dm}{d\pi}\right)dm, 
 \end{align} 
where $\frac{dm}{d\pi}$ is the Radon-Nikodym derivative of $m$ with respect to $\pi$~\cite{yeh---14}.

Rather than giving a definition of the conditional mutual information $I(x;y|z)$ in terms of conditional distributions, it suffices for our purposes to define it through the chain rule of mutual information\cite[Corollaries Corollary~2.5.1 and~7.14]{gray--11}:
\begin{align}
 I(x;y|z)\eq I(x;y,z)-I(x;z)=I(x,z;y)-I(z;y).
\end{align}
If $\rvay\in\Asp$ is a discrete-valued random variable with alphabet $\Asp$, then
\begin{align}\label{eq:ictohc}
 I(\rvax;\rvay|\rvaz)=H(\rvay|\rvaz)-H(\rvay|\rvax,\rvaz),
\end{align}
where
\begin{align}
 H(\rvay|\rvaz)\eq -\Expe{\Sumover{u\in\Asp}\Pr\set{\rvay=u|\rvaz}\log_{2}\left(\Pr\set{\rvay=u|\rvaz}\right)}
\end{align}
is the conditional entropy of $\rvay$ given $\rvaz$, 
$\Pr\set{\rvay=u|\rvaz}$ is the probability that $\rvay=u$ given $\rvaz$ and the expectation $\Expe{}$ is with respect to $\rvaz$.

\begin{defn}[Directed Information with Forward Delay]
\textit{In this paper, the directed information from $\rvax^{k}$ to $\rvay^{k}$ through a forward channel with a non-negative time varying delay of $d_{xy}(i)$ samples is defined as
\begin{align}
 I(\rvax^{k}\to \rvay^{k}) \eq \Sumfromto{i=0}{k}I(\rvay(i);\rvax^{i-d_{xy}(i)}|\rvay^{i-1}).
\end{align}}
\end{defn}
For a zero-delay forward channel, the latter definition coincides with Massey's~\cite{massey90}.

Likewise, we adapt the definition of causally-conditioned directed information to the definition
\begin{align}\label{eq:dimc}
 I(\rvax^{k}\to \rvay^{k}\parallel\rvae^{k} ) 
 \eq 
 \Sumfromto{i=0}{k}I(\rvay(i);\rvax^{i-d_{xy}(i)}|\rvay^{i-1},\rvae^{i}).
\end{align}
where, as before, $d_{xy}(i)$ is the delay from $\rvax$ to $\rvay(i)$.

\section{Theorem Statement and the Flaw in its Proof}\label{ssec:flaw}

\begin{figure}[htbp]
\centering
\input{NCS_from_framweork_paper.pstex_t}
\caption{General source-coding scheme for the networked control of an LTI plant $\mathscr{F}$ considered in~\cite{silder11}.}
\label{fig:e-d}
\end{figure}
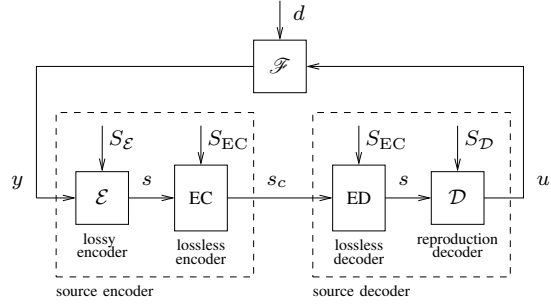
The networked control system considered in~\cite{silder11}, depicted in Fig.~\ref{fig:e-d}, considers lossy encoders
$\enc$, reproduction decoders $\dec$ and EC-ED pairs that are
causal, operating without delay. 
It is also considered the possible availability of side information 
sequences  
$S_{\enc}(k)$ at the encoder and 
$S_{\dec}(k)$ at the decoder.  
Such side information is contained in suitably defined sets $\mathcal{S}_{\dec}(k)$ and
$\mathcal{S}_{\enc}(k)$, where $S_{\enc}(k) \in
\mathcal{S}_{\enc}(k)$ and $S_{\dec}(k) \in
\mathcal{S}_{\dec}(k)$. 
Next,~\cite{silder11} defines the set $\mathcal{S}_{\EC}(k)\triangleq \mathcal{S}_{\enc}(k) \cap
\mathcal{S}_{\dec}(k)$, which contains the common side
information that becomes available at both the encoder and
decoder sides at instant $k$.
More explicitly, what the authors of~\cite{silder11} meant  was
\begin{align}
 \mathcal{S}_{\enc}(k)&\eq \Set{S_{\enc}(k) ,S_{\EC}(k)  }, \fspace k=0,1,\ldots
\\  
\mathcal{S}_{\dec}(k)&\eq \Set{S_{\dec}(k) ,S_{\EC}(k)  }, \fspace k=0,1,\ldots\label{eq:blg}
\end{align}
 
The output $y$ of the plant $\pla$ satisfies
\begin{subequations}\label{eq:gen}
\begin{align}\label{eq:feedback}
y(k)=\mathscr{F}_k(u^{k-1}, d^k, x_o),
\end{align}
where $\mathscr{F}_k$ is a (possibly time varying) deterministic mapping,
$x_o\in\Rl^{n_o}$ is the initial state of $\mathscr{F}$, $d$, with $d(k)\in\Rl^{n_d}$, is an
exogenous random process 
and $u$ is the plant input.
In~\cite{silder11}, $u$ and $y$ were real valued scalars, but for its Theorem~4.1 and for this note, their samples may take their outcomes from arbitrary sets.
The other sequences are generated as
\begin{align}\label{eq:enc}
s(k)&= \enc_{k}(y^k, S_{\enc}^k) \in \Asp_{s},
\\
s_c(k) &= \mathcal{H}_k ( s^{k}, S_{\EC}^k ) \in\Asp(k),
\\
s(k) &= \mathcal{H}_k^{-1} (s_c^k, S_{\EC}^k ),
\\
u(k)&= \dec_{k}(s^k, S_{\dec}^k),
\end{align}
\end{subequations}
where
$\mathcal{A}_s$ is a fixed countable
set and
 $\mathcal{A}(k)$ is a countable set of prefix-free binary
words.

The expected length $R(k)$ of any binary description $s_c(k)$
of the lossy encoder output symbol $s(k)$ satisfies (see
\cite[Chapter~5]{covtho06},~\cite{shanno48} and also
~\cite{lintar97})
\begin{align}\label{eq:fundamental-R}
H(s(k) | s^{k-1}, S_{\EC}^k ) \leq  R(k).
\end{align}

From these definitions and results,~\cite[Theorem~4.1]{silder11} can be stated as follows
\begin{thm}[{\cite[Theorem~4.1]{silder11}}]\label{thm:silder4.1} 
 Consider a source-coding scheme inside a feedback loop, as depicted in Fig.~\ref{fig:e-d}, where the sequences are generated via~\eqref{eq:gen}.
If 
$\Ssp_{\dec}^{k}\Perp (x_{o},d^{k})$,
and $\Dsp$ is conditionally invertible, 
then, for all $k\in\Nl_{0}$,
$\sumfromto{i=0}{k}R(i) \geq I(y^{k}\to u^{k})$.
Thus,
\begin{align}
 \lim_{k\to\infty}\frac{1}{k+1}
 \sumfromto{i=0}{k}R(i) \geq 
  \lim_{k\to\infty}\frac{1}{k+1}
I(y^{k}\to u^{k})
\end{align}
provided both limits exist.
\end{thm}
This is a key result, because, combined with~\cite[eq.~9]{silder11}, 
it yields
\begin{align}\label{eq:sdw}
  \frac{1}{k+1} I(y^{k}\to u^{k})
  \leq\frac{1}{k+1} \Sumfromto{i=0}{k} R(i)\leq \frac{1}{k+1} I(y^{k}\to u^{k})+1 \fspace\text{[bits/sample]}, \fspace k=0,1,2,\ldots.
\end{align}
This result highlights the operational meaning of the directed information as a lower bound (tight to within one bit) to the data rate of any given source code in a closed-loop system.
This fact has been a crucial ingredient in characterizing the best 
rate-performance achievable in Gaussian linear quadratic networked 
control~\cite{silder16,tanesf18}.

The proof of Theorem 4.1 in~\cite{silder11} is invalid since it relies upon~\cite[Lemma 4.2]{silder11}, whose first claim does not hold.
 The latter claim was that the Markov chain 
$
S_{\dec}^{i} 
\longleftrightarrow 
u^{i-1} 
\longleftrightarrow 
y^{i}
$
held.
The flawed reasoning in the proof of this result is the following:
\begin{quotation}
``Given $u^{i-i}$, it follows from~\eqref{eq:feedback} that there exists a deterministic mapping $T_{i}$ such that $y^{i}=T_{i}(d^{i},x_{o})$.
Since $S_{\dec}^{i}\Perp (d^{i},x_{o})$, it immediately follows that $y^{i}$ and $S_{\dec}^{i}$ are independent upon knowledge of $u^{i-i}$''.
\end{quotation}
The problem with this  argument is that it neglects the relationship between $S_{\dec}^{i}$ and $y^{i}$ through $s^{i}$.
Specifically, given $u^{i-1}$, one may gain information about $s^{i-1}$ from $S_{\dec}^{i}$.
 This, in turn, may give information about $y^{i}$.
 More generally, two independent random variables may cease to be independent given a third one.

\section{A More General Theorem}\label{ssec:flawfix}
\begin{figure}[htbp]
\centering
\begin{tikzpicture}[node distance=9mm]
 \def\la{3.5mm}
\node (b1) [bloq] {$\Ssp_{1}$};
\node (b2) [bloq, right =of b1] {$\Ssp_{2}$};
\node (b3) [bloq, below = of b2, yshift=3mm] {$\Ssp_{3}$};
\node (b4) [bloq, left =of b3] {$\Ssp_{4}$};

\path (b1) edge[->] node[pos=0.5, above] {$\rvae$} (b2);
\path (b3) edge[->] node[pos=0.5, below] {$\rvay$} (b4);
\draw[->] (b2.east) --++ (5mm,0) |-($(b3.east) $) node[pos=.75,below] {$\rvax$};
\draw[->] (b4.west) --++(-.5,0) |-(b1.west) node[ pos=.75, above] {$\rvau$};

\draw[<-] (b1.north)--++(0,\la) node[above] {$\rvar$};
\draw[<-] (b2.north)--++(0,\la) node[above] {$\rvap$};
\draw[<-] (b3.south)--++(0,-\la) node[below] {$\rvas$} coordinate (ps);
\draw[<-] (b4.south)--++(0,-\la) node[below] {$\rvaq$} coordinate (pq);
\end{tikzpicture}
\caption{The general system considered by Theorem~\ref{thm:derostthm2}}.
\label{fig:diagramas}
\end{figure}
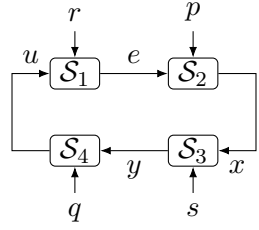

In this section we dorive a proof for Theorem~\ref{thm:silder4.1} which does not require a conditionally invertible decoder.
This result builds upon a recent theorem by the authors that applies to the general feedback system shown in Fig.~\ref{fig:diagramas}.
In this diagram, the blocks $\Ssp_{1},\ldots, \Ssp_{4}$ represent possibly non-linear and time-varying causal discrete-time systems such that the total delay of the loop is at least one sample. 
In the same figure, $\rvar,\rvap,\rvas,\rvaq$ are exogenous random signals (scalars, vectors or sequences), which could represent, for example, any combination of disturbances, noises, random initial states or side informations.
Note that any of these exogenous signals, in combination with its corresponding deterministic mapping $\Ssp_{i}$, can yield any desired stochastic causal mapping (for example, a noisy communication channel, a zero-delay source coder or decoder, or a causal dynamic system with disturbances and a random initial state).

\begin{thm}[Full Conditional Closed-Loop Directed Data-Processing Inequality from~{\cite[Thm.~2]{derost21a}}]\label{thm:derostthm2}
Consider  the system shown in Fig.~\ref{fig:diagramas}.

If $(\rvaq,\rvas)\Perp (\rvar,\rvap)$ and%
 \footnote{The Markov chain notation $\rvaa\leftrightarrow \rvab \leftrightarrow \rvac$ means ``$\rvaa$ and $\rvac$ are independent when $\rvab$ is given''.}
 $\rvaq_{i+1}^{k}\leftrightarrow \rvaq^{i} \leftrightarrow \rvas^{i}$ for $i=0,1,\ldots,k-1$, then
\begin{align}\label{eq:silderok}
 I(\rvax^{k}\to\rvay^{k}\Vert\rvaq^{k})\geq I(\rvae^{k}\to\rvau^{k}).
\end{align}
\finenunciado
\end{thm}

Now we can prove the following more general version of Theorem~\ref{thm:silder4.1}.
\begin{thm}\label{thm:new}
Consider a source-coding scheme inside a feedback loop, as depicted in Fig.~\ref{fig:e-d}, where the sequences are generated via~\eqref{eq:gen}.
If 
$(\Ssp_{\enc}^{k},\Ssp_{\dec}^{k})\Perp (x_{o},d^{k})$,
and
\begin{align}\label{eq:mcsi}
(\Ssp_{\dec}(i+1),\Ssp_{\dec}(i+2),\ldots,\Ssp_{\dec}(k))\leftrightarrow\Ssp_{\dec}^{i}\leftrightarrow\Ssp_{\enc}^{i}
, \fspace i=0,1,\ldots,k-1,
\end{align}
then, for all $k\in\Nl_{0}$,
$\sumfromto{i=0}{k}R(i) \geq I(y^{k}\to u^{k})$.
Thus,
\begin{align}
 \limsup_{k\to\infty}\frac{1}{k+1}\sumfromto{i=0}{k}R(i) &\geq \limsup_{k\to\infty}\frac{1}{k+1}I(y^{k}\to u^{k})
\\  
\liminf_{k\to\infty}\frac{1}{k+1}\sumfromto{i=0}{k}R(i) &\geq \liminf_{k\to\infty}\frac{1}{k+1}I(y^{k}\to u^{k})
\end{align}
\end{thm}

\begin{proof}
 We have that 
\begin{align}
 \Sumfromto{i=0}{k}
R(i)
&
\overset{(a)}{\geq}
 \Sumfromto{i=0}{k}
H( s(i)|  s^{i-1},\Ssp_{\dec}^{k})
\\&
\overset{(b)}{\geq}
 \Sumfromto{i=0}{k}\left(
H( s(i)|  s^{i-1},\Ssp_{\dec}^{i})
-
H( s(i)|  s^{i-1},\Ssp_{\dec}^{i}, y^{i})\right)
\\&
\overset{\eqref{eq:ictohc}}{=}
 \Sumfromto{i=0}{k}
I( s(i); y^{i}|  s^{i-1},\Ssp_{\dec}^{i})
\\&
\overset{\eqref{eq:dimc}}{=}
I( y^{k} \to  s^{k}\Vert\Ssp_{\dec}^{k})
\\&
\geq
I( y^{k} \to \rvau^{k})
\end{align}
where 
$(a)$ follows from~\eqref{eq:fundamental-R} and~\eqref{eq:blg},~$(b)$ is from the non-negativity of entropy and
the last inequality follows directly from Theorem~\ref{thm:derostthm2}.
\end{proof}

\begin{rem}
 The additional independency assumptions made in Theorem~\ref{thm:new} with respect to Theorem~\ref{thm:silder4.1} are not a problem for the other results in~\cite{silder11} that rely upon Theorem~4.1.
 This is so because they are met by the only side information signals considered in~\cite{silder11}:
 identical dither sequences available to the encoder and decoder.
 \finenunciado
\end{rem}

\section{Conclusions}
We have revealed a flaw in the proof of Theorem~4.1 in~\cite{silder11} and derived a proof for it without requiring an invertible decoder.
Thus, we have extended the validity of~\cite[Theorem~4.1]{silder11} in such a way that all the other results in~\cite{silder11} which depend on that theorem continue to hold.

% %--------------------------------------------------------------------
% % THE END
% %--------------------------------------------------------------------
%-----------------------------------------------------------------

\end{document}

%% file: NCS_from_framweork_paper.pstex_t
\begin{picture}(0,0)%
\includegraphics{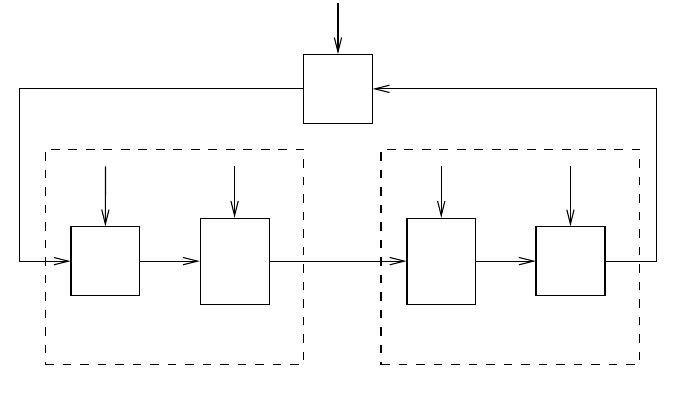}%
\end{picture}%
\setlength{\unitlength}{2763sp}%
\begingroup\makeatletter\ifx\SetFigFont\undefined%
\gdef\SetFigFont#1#2#3#4#5{%
  \reset@font\fontsize{#1}{#2pt}%
  \fontfamily{#3}\fontseries{#4}\fontshape{#5}%
  \selectfont}%
\fi\endgroup%
\begin{picture}(4636,2684)(872,-3605)
\put(3958,-2232){\makebox(0,0)[lb]{\smash{{\SetFigFont{8}{9.6}{\rmdefault}{\mddefault}{\updefault}{\color[rgb]{0,0,0}$S_{\operatorname{EC}}$}%
}}}}
\put(2540,-2232){\makebox(0,0)[lb]{\smash{{\SetFigFont{8}{9.6}{\rmdefault}{\mddefault}{\updefault}{\color[rgb]{0,0,0}$S_{\operatorname{EC}}$}%
}}}}
\put(3899,-2764){\makebox(0,0)[b]{\smash{{\SetFigFont{7}{8.4}{\rmdefault}{\mddefault}{\updefault}{\color[rgb]{0,0,0}ED}%
}}}}
\put(4253,-2586){\makebox(0,0)[lb]{\smash{{\SetFigFont{8}{9.6}{\rmdefault}{\mddefault}{\updefault}{\color[rgb]{0,0,0}$s$}%
}}}}
\put(3899,-3177){\makebox(0,0)[b]{\smash{{\SetFigFont{6}{7.2}{\rmdefault}{\mddefault}{\updefault}{\color[rgb]{0,0,0}lossless}%
}}}}
\put(3899,-3295){\makebox(0,0)[b]{\smash{{\SetFigFont{6}{7.2}{\rmdefault}{\mddefault}{\updefault}{\color[rgb]{0,0,0}decoder}%
}}}}
\put(1595,-2764){\makebox(0,0)[b]{\smash{{\SetFigFont{8}{9.6}{\rmdefault}{\mddefault}{\updefault}{\color[rgb]{0,0,0}$\enc$}%
}}}}
\put(1950,-2586){\makebox(0,0)[lb]{\smash{{\SetFigFont{8}{9.6}{\rmdefault}{\mddefault}{\updefault}{\color[rgb]{0,0,0}$s$}%
}}}}
\put(1596,-3118){\makebox(0,0)[b]{\smash{{\SetFigFont{6}{7.2}{\rmdefault}{\mddefault}{\updefault}{\color[rgb]{0,0,0}lossy }%
}}}}
\put(1596,-3236){\makebox(0,0)[b]{\smash{{\SetFigFont{6}{7.2}{\rmdefault}{\mddefault}{\updefault}{\color[rgb]{0,0,0}encoder}%
}}}}
\put(3072,-2586){\makebox(0,0)[lb]{\smash{{\SetFigFont{8}{9.6}{\rmdefault}{\mddefault}{\updefault}{\color[rgb]{0,0,0}$s_c$}%
}}}}
\put(3190,-1582){\makebox(0,0)[b]{\smash{{\SetFigFont{8}{9.6}{\rmdefault}{\mddefault}{\updefault}{\color[rgb]{0,0,0}$\mathscr{F}$}%
}}}}
\put(3308,-1110){\makebox(0,0)[lb]{\smash{{\SetFigFont{8}{9.6}{\rmdefault}{\mddefault}{\updefault}{\color[rgb]{0,0,0}$d$}%
}}}}
\put(5493,-2586){\makebox(0,0)[lb]{\smash{{\SetFigFont{8}{9.6}{\rmdefault}{\mddefault}{\updefault}{\color[rgb]{0,0,0}$u$}%
}}}}
\put(887,-2586){\makebox(0,0)[rb]{\smash{{\SetFigFont{8}{9.6}{\rmdefault}{\mddefault}{\updefault}{\color[rgb]{0,0,0}$y$}%
}}}}
\put(1655,-2232){\makebox(0,0)[lb]{\smash{{\SetFigFont{8}{9.6}{\rmdefault}{\mddefault}{\updefault}{\color[rgb]{0,0,0}$S_{\enc}$}%
}}}}
\put(1182,-3590){\makebox(0,0)[lb]{\smash{{\SetFigFont{6}{7.2}{\rmdefault}{\mddefault}{\updefault}{\color[rgb]{0,0,0}source encoder}%
}}}}
\put(3485,-3590){\makebox(0,0)[lb]{\smash{{\SetFigFont{6}{7.2}{\rmdefault}{\mddefault}{\updefault}{\color[rgb]{0,0,0}source decoder}%
}}}}
\put(2481,-2763){\makebox(0,0)[b]{\smash{{\SetFigFont{7}{8.4}{\rmdefault}{\mddefault}{\updefault}{\color[rgb]{0,0,0}EC}%
}}}}
\put(2481,-3177){\makebox(0,0)[b]{\smash{{\SetFigFont{6}{7.2}{\rmdefault}{\mddefault}{\updefault}{\color[rgb]{0,0,0}lossless}%
}}}}
\put(2481,-3295){\makebox(0,0)[b]{\smash{{\SetFigFont{6}{7.2}{\rmdefault}{\mddefault}{\updefault}{\color[rgb]{0,0,0}encoder}%
}}}}
\put(4784,-2763){\makebox(0,0)[b]{\smash{{\SetFigFont{8}{9.6}{\rmdefault}{\mddefault}{\updefault}{\color[rgb]{0,0,0}$\dec$}%
}}}}
\put(4784,-3118){\makebox(0,0)[b]{\smash{{\SetFigFont{6}{7.2}{\rmdefault}{\mddefault}{\updefault}{\color[rgb]{0,0,0}reproduction}%
}}}}
\put(4784,-3236){\makebox(0,0)[b]{\smash{{\SetFigFont{6}{7.2}{\rmdefault}{\mddefault}{\updefault}{\color[rgb]{0,0,0}decoder}%
}}}}
\put(4844,-2232){\makebox(0,0)[lb]{\smash{{\SetFigFont{8}{9.6}{\rmdefault}{\mddefault}{\updefault}{\color[rgb]{0,0,0}$S_{\dec}$}%
}}}}
\end{picture}%